\documentclass[12pt,leqno]{article}
\usepackage{amsfonts}
\usepackage{ccaption,mathrsfs,amsmath}
\usepackage{amsmath, amscd, amsfonts, amsthm, tikz, times}

\newcounter{conjecture}
\newcounter{remark}
\newcounter{corollary}

\newcommand{\eqnsection}{
    \renewcommand{\theequation}{\thesection.\arabic{equation}}
    \makeatletter
    \csname @addtoreset\endcsname{equation}{section}
    \makeatother}
\newtheorem{theorem}{Theorem}
\newtheorem{lemma}{Lemma}
\newtheorem{deff}{Definition}
\newtheorem{prop}{Proposition}
\newtheorem{cor}{Corollary}

\newcommand{\ind}{\hspace{.25in}}
\newcommand{\dd}{\delta}
\newcommand{\DD}{\Delta}
\newcommand{\subs}{\subseteq}
\newcommand{\lar}{\longrightarrow}
\newcommand{\eps}{\varepsilon}

\newcommand{\reals}{R}
\newcommand{\CC}{\mathbb{C}}

\newcommand{\lll}{\label}

\def \be{\begin{equation}}
\def \ee{\end{equation}}
\def \bt{\begin{theorem}}
\def \et{\end{theorem}}
\def \bea{\begin{eqnarray}}
\def \eea{\end{eqnarray}}
\def \bas{\begin{eqnarray*}}
\def \eas{\end{eqnarray*}}

\newcommand {\rrr}[1]{(\ref{#1})}

\def \la{\lambda}

\def \ff{\infty}

\def \DD{\mathbb{D}}

\def \TT{{\cal T}}

\def \({\left(}
\def \){\right)}

\def \nn{\nonumber}

\def \vski{\vspace{12pt}}

\def \bc{\begin{center} }
\def \ec{\end{center} }
\def \bs{\begin{slide} }
\def \es{\end{slide} }

\def\square{{\vcenter{\vbox{\hrule height.3pt
         \hbox{\vrule width.3pt height5pt \kern5pt
            \vrule width.3pt}
         \hrule height.3pt}}}}
\def\qed{{\hfill $\Box$ \bigskip}}

\eqnsection
\begin{document}

\title{The rate of convergence of the hyperbolic density on sequences of domains.}



\author{
\begin{tabular}{cc}
\textit{Nikola Lakic} & \textit{Greg Markowsky} \\
\end{tabular}}

\bibliographystyle{amsplain}
\maketitle \eqnsection \setlength{\unitlength}{2mm}

\begin{abstract}
It is known that if a sequence of domains $U_n$ converges to a domain $U$ in the Carath\'{e}odory sense then the hyperbolic densities on $U_n$ converge to the hyperbolic density on $U$. In this paper, we study the rate of convergence of the hyperbolic density under a slightly different mode of convergence. In doing so, we are led to consider two other densities on domains, the Teichm\"{u}ller density and the three-point density. We obtain several results which give rates of convergence in various scenarios.
\end{abstract}

\let\thefootnote\relax

\let\thefootnote\relax\footnote{The first author was partially supported by NSF grant 0700052 and PSC-CUNY grant 63114-00-41. The second author was supported by the Priority Research Centers Program through the National Research Foundation of Korea (NRF) funded by the Ministry of Education, Science and Technology (MEST) (No. 2010-0029638) as well as Australian Research Council Grant DP0988483.}

\section{Introduction}

\ind The well-known hyperbolic, or Poincar\'{e}, density on plane domains has proved to be of great importance and utility in complex analysis and geometry. It is known that if $U$ be a hyperbolic domain and $U_n$ a sequence of domains
in $\CC$ which converge to $U$ in a reasonable way, then the hyperbolic density on $U_n$ converges locally uniformly to that on $U$. However, for applications it may be useful to be able to say something about the rate of this convergence. This is the question, originally posed by A. Douady, to which this paper is devoted.

\ind Below, all distances are measured to the Euclidean metric unless otherwise specified. The Hausdorff distance $H(A,B)$ between any
two sets $A$ and $B$ is defined to be

\be \label{}
H(A,B):= \inf_{r>0} \{A \subs N_r(B) \mbox{ and } B \subs N_r(A) \},
\ee

where $N_r(A)$ is the set of all points whose
distance from $A$ is less than $r$. In contrast, we define

\be \label{}
d(A,B):= \inf_{z\in A,w \in B} \{|z-w|\},
\ee

For example, if we let $A=\{1/2\}, B = \{|z|=1\}$, and $C=\{|z-2|=1\}$, then $H(A,B)=3/2, H(A,C)=5/2, H(B,C)=2, d(A,B)=1/2, d(A,C)=1/2, d(B,C)=0$. For singleton sets we will commonly use the shorthands $H(x,A), d(x,A)$ instead of $H(\{x\},A), d(\{x\},A)$. The relation between $d$ and $H$ is that

\begin{equation} \label{}
H(A,B) = \max\Big(\sup_{z \in A} d(z,B) , \sup_{z \in B} d(z,A) \Big)
\end{equation}

The following is the definition of convergence of domains which will be useful to us.

\begin{deff} \lll{d1} Let $U$ be a domain and $U_n$ a sequence of domains
in $\CC$. We say that $U_n$ converges in boundary to $U$ if

\begin{itemize}
\item[a)]{$H(\dd U, \dd U_n) \lar 0$}

\item[b)]{There exists $z_0 \in U$ such that $z_0 \in U_n$ for all $n$.}
\end{itemize}

\end{deff}

{\bf Remark:} Condition $(b)$ is necessary to
eliminate such situations as $U=\DD$, $U_n = \{1-1/n < |z| <
1\}$.

\vski

A related concept is {\it Carath\'{e}odory convergence}. Consider a sequence of domains $\{U_n\}$ containing the fixed point $z_0$. The {\it kernel} $U$ of this sequence is defined to be the maximal domain containing $z_0$ such that every compact $K \subseteq U$ is contained in $U_n$ for sufficiently large $n$. If no such $U$ exists, then let $U = \{z_0\}$.  We say that $U_n$ converges to its kernel $U$ in the Carath\'{e}odory sense if $U$ is also the kernel of every subsequence of $\{U_n\}$ (see \cite{dur} or any of a number of other books on univalent functions or geometric function theory). The following lemma shows the connection between the two modes of convergence.

\begin{lemma} If $U_n \lar U$ in boundary then $U_n \lar U$ in the Carath\'{e}odory sense.
\end{lemma}

{\bf Proof:} Let $z_0$ be as in Definition \ref{d1}, and suppose $K \subseteq U$ is compact. We may assume $K$ is connected and contains $z_0$, as it will
always be contained in a compact, connected set within $U$ containing $z$. By
compactness, $K$ is a positive distance from $\dd U$. Condition $(a)$
then implies that, for sufficiently large $n$, $K \bigcap \dd U_n =
\emptyset$, and thus $K$ lies in $(\dd U_n)^c$. Since $K$ is
connected, it lies within a connected component of $(\dd U_n)^c$. As
$z \in K$ and $z \in U_n$ for large $n$, it follows that the
connected component in question is $U_n$ itself. This shows that the kernel of $\{U_n\}$ is contained in $U$. However, convergence in boundary easily implies

\begin{equation} \label{}
\bigcap_{N=1}^\ff \bigcup_{n \geq N} U_n \subseteq \bar{U} .
\end{equation}

Thus, $U$ is the kernel of $\{U_n\}$. It is clear that any subsequence of $\{U_n\}$ also converges in boundary to $U$, and thus by the same argument has $U$ as its kernel. We see that $U_n \lar U$ in the Carath\'{e}odory sense. \qed

The converse to this lemma is false, as the domains $U_n=\DD \bigcup \{z \in \CC: \arg(z) \in (0,\frac{1}{n})\}$ converge to $U=\DD$ in the Carath\'{e}odory sense but not in boundary. Since we are interested in the behavior of the hyperbolic density we assume henceforth that $U$ has at least 3 boundary points in the Riemann sphere, which implies that the hyperbolic density for $U$ exists. This clearly implies that if $U_n \lar U$ in boundary, then $U_n$ is hyperbolic for sufficiently large $n$. Let $\rho_A(z)$ denote the hyperbolic density of any hyperbolic domain $A$ at the point $z$, normalized to have curvature -4 (a normalization with curvature -1 is also common). Then $\rho_A(z) = \frac{1}{|\pi '(0)|}$, where $\pi$ is a holomorphic covering map from $\DD$ to $A$ with $\pi(0)=z$. The following lemma should be considered known, though it may not yet have been stated in this form.

\begin{lemma} \label{old} If $U_n \lar U$ in boundary, then $\rho_{U_n}(z) \lar \rho_U (z)$ locally uniformly.
\end{lemma}

{\bf Proof:} This is immediate from Theorem 1 of \cite{hej}, which shows under the more general condition of Carath\'{e}odory convergence that the covering maps from $\DD$ to $U_n$ converge locally uniformly to that of $U$. \qed

{\bf Remark:} For a given $z$, $\rho_{U_n}(z)$ may only be defined for large $n$, and the
statement of the theorem should be interpreted accordingly.

\section{The Teichm\"{u}ller density} \lll{tm}

\ind In order to say something about the rate of convergence of $\rho_{U_n}(z)$ to $\rho_{U}(z)$, we introduce a different conformally invariant density. Denote the Teichm\"{u}ller density on a domain $A$ by $\la_A (z)$, defined as

\be \label{}
\la_A (z) := \inf_{V \in \TT} ||\bar{\dd} V||_\ff
\ee

where $||\cdot||_\ff$ denotes the $L^\ff$ norm, $\bar{\dd}V = \frac{dV}{d\bar{w}} = \frac{1}{2}(\frac{dV}{dx} + i\frac{dV}{dy})$, and $\TT$ is the family of all complex-valued functions $V(w)$ on $A$ with $V(z) = 1$, which vanish on $\dd A$, and with distributional derivatives. Since $\bar{\dd}V$ is bounded, the integral $\int_A \bar{\dd}V \varphi$ converges for every integrable function $\varphi$ on $A.$  By normalizing, we may temporarily assume that 0 and 1 are in the complement of $A.$   Cauchy's formula shows $V(z)=\int_A \bar{\dd}V \varphi_z$ where $\varphi_z(w)=-\frac{z(z-1)}{\pi w(w-1)(w-z)}.$   The linear span of functions of the form $\varphi_z$ where $z$ is in the complement of $A - \{p\}$ is dense in the space of all integrable and holomorphic functions on  $A - \{p\}.$  Thus, the chain rule applied to $W(w)=V(f(w))\frac{f'(z)}{f'(w)}$
shows that $\la$ is conformally invariant, in the sense that if $f$ is a conformal map we have $\la_A(z)=\la_{f(A)}(f(z))|f'(z)|$. The Teichm\"{u}ller density has still at this point a relatively brief history. It was originally defined in \cite{fnart}, where it was used in the study of uniformly thick and uniformly perfect domains. It was also proved there that the Teichm\"{u}ller and hyperbolic densities are equivalent, in the sense that $\frac{1}{2}\rho_A \leq \la_A \leq \rho_A$ for any domain $A$. The transitivity of the automorphism group of $\DD$ together with the conformal invariance of $\la_A$ shows that $\la_A$ and $\rho_A$ coincide, up to a multiplicative constant, on simply connected domains. It is possible to calculate this constant, and it turns out that $\frac{1}{2}\rho_A = \la_A$ in this case. On the other hand, if $A$ is the thrice punctured sphere then it is known that $\rho_A = \la_A$ (see \cite{egn}). It was shown in \cite{yu} that $\la_A$ is continuous for any domain $A$, and is the infinitesimal form of a previously known metric defined in terms of Teichm\"{u}ller shift mappings(see \cite{kra}).

\ind We now give an intuitive explanation of the Teichm\"{u}ller density. Holomorphic functions are functions which satisfy $\bar{\dd}f = 0$, so in essence $\la_A$ is measuring how nearly holomorphic a function can be while attaining the prescribed values at $z$ and on the boundary. The definition of the Teichm\"{u}ller density arises most naturally in the context of holomorphic motions. Given any closed set $E$ in the extended complex plane $\hat{\CC}$, a holomorphic motion $h_t(w)$ is a function from $\DD \times E \lar \hat{\CC}$ which satisfies the following properties.

\begin{itemize}
\item[i)]{$t \lar h_t(w)$  is holomorphic on $\DD$ for every fixed $w \in E$.}
\item[ii)]{$w \lar h_t(w)$ is injective on $E$ for every fixed $t \in \DD$.}
\item[iii)]{$f_0(w) = w$ for every $z \in E$.}
\end{itemize}

Note the lack of any sort of continuity assumption on $h$ as a function of $w$. For this reason, it may be helpful to suppose first that $E$ is a finite set and that $h$ effects a simultaneous motion of the points in $E$. This motion has a complex time variable $t$, and the points in $E$ are not allowed to collide at any time. However, it is a remarkable fact that the holomorphicity in $t$ forces $h$ to satisfy strong continuity conditions in $w$ if $E$ is an infinite set. In fact, much more is true. The following statement, commonly referred to as the $\la-lemma$, was first proved by Slodkowski, although a weaker version had been proved earlier by Sullivan and Thurston (see \cite{fyun} for a complete account).

\begin{lemma}
Suppose $E$ is a closed set, and $h_t(w): \DD \times E \lar \hat{\CC}$ is holomorphic motion. Then there is a holomorphic motion $\tilde{h}_t(w): \DD \times \hat{\CC} \lar \hat{\CC}$ such that $\tilde{h}$ agrees with $h$ on $\DD \times E$. For fixed $t$, the function $w \lar \tilde{h}_t(w)$ is a quasiconformal homeomorphism from $\hat{\CC}$ to $\hat{\CC}$.
\end{lemma}

Note that nothing is said about uniqueness, and in general there will be many possible extensions of a given $h$. Suppose $E=U^c \bigcup \{z\}$, and suppose further that we have a holomorphic motion $h_t(w)$ such that $\frac{d}{dt} h_t(w)\Big|_{t=0} = 0$ for all $w \in U^c$ and $\frac{d}{dt} h_t(z)\Big|_{t=0} = 1$. Let $\tilde{h}$ be an extension of this holomorphic motion to all of $\hat{\CC}$ guaranteed by Slodkowski's Theorem, and let us consider $\phi:=\frac{d^2}{d\bar{w}dt} \tilde{h}_t(w)\Big|_{t=0}$. Taking the $t$ derivative first and setting $V(w)=\frac{d}{dt}\tilde{h}_t(w)\Big|_{t=0}$, we obtain $\phi = \frac{dV}{d\bar{w}}$. On the other hand, $\tilde{h}_t(w)$ is quasiconformal in $w$ and thus has a Beltrami coefficient $\mu_t(z)$ such that $\frac{d\tilde{h}}{d\bar{w}} = \mu_t \frac{d\tilde{h}}{dw}$. We then have

\be \label{}
\frac{d\mu_t}{dt} = \frac{\frac{d\tilde{h}}{dw}\frac{d^2\tilde{h}}{dtd\bar{w}}-\frac{d\tilde{h}}{d\bar{w}}\frac{d^2\tilde{h}}{dtdw}}{\Big(\frac{d\tilde{h}}{dw}\Big)^2}
\ee

Note that $\frac{d\tilde{h}}{d\bar{w}}(w)\Big|_{t=0}=0$ and $\frac{d\tilde{h}}{dw}(w)\Big|_{t=0}=1$ for all $w$, so that $\phi = \frac{d^2}{d\bar{w}dt} \tilde{h}_t(w)\Big|_{t=0} = \frac{d\mu_t}{dt}\Big|_{t=0}$. Comparing our two expressions for $\phi$, we see $\frac{dV}{d\bar{w}}(z) = \frac{d\mu_t}{dt}(z)\Big|_{t=0}$. Thus, the Teichm\"{u}ller density measures the minimal rate of change of the Beltrami coefficient of all quasiconformal homeomorphisms associated to a holomorphic motion which fixes the boundary of $U$ and moves $z$ with unit velocity at time $0$. It stands to reason that points closer to the boundary, in whatever sense, require a more violent holomorphic motion in order to move (while the boundary remains fixed) than those farther away. This results in a larger value of $\la_U$, as is the case with the hyperbolic metric $\rho_U$.

\ind This density is well suited to our problem concerning the rate of convergence of densities on domains, as we can use known moduli of continuity on vector fields associated to holomorphic motions to our advantage. Suppose that a sequence of domains $U_n$ converges in boundary to a domain $U$. The following theorem shows that the Teichm\"{u}ller metric on $U_n$ converges to the Teichm\"{u}ller metric on $U$ with speed $\log(\log(\frac{1}{H(\dd U_n,\dd U)}))$.

\begin{theorem} \lll{big}
Suppose that $U$ is a bounded domain, and that $K$ is a compact
subset of $U$. Then there is a constant $C$ such that if $W$ is a
domain with $H(\dd U,\dd W) \leq \eps$ then $|\la_U (w) - \la_W (w)|
\leq C (\log(\log(1/\eps)))^{-1}$ for all $w \in K$, provided that
$\eps$ is sufficiently small. The constant $C$ depends on $d(K,\dd U)$ and the diameter of $U$.
\end{theorem}

{\bf Proof:} Let $\eps < \eps_o<e^{-1}$ where $\eps_o$ is sufficiently
small so that the condition $d(\dd U,\dd W) \leq \eps_o$ forces $W$
to contain $K$. Let $d_U (z) = d(z,\dd U)$. Given $w \in K$,
let $V$ denote a differentiable vector field on $U$ with $V(w) = 1$,
$V(z) = 0$ for $z \in U^c$, and $||\bar{\dd}V||_{\infty} = \la_U
(w)$. We will use $V$ to construct a vector field on $W$ which will
give us a bound for $\la_W (w)$. The problem, of course, is that $V(z)$
is not necessarily 0 on $W^c$, so we must alter it at the boundary
in a way that doesn't affect the norm of the $\bar{\dd}$ derivative
very much. We can assume, after multiplying the entire picture by a
constant if necessary, that $d(\dd U,K) \geq 1/2$. Let

\be
j(x)= \left \{
\begin{array}{ll}
0, & \qquad \mbox{if { }} 0 \leq x \leq \eps, \\
(\log(\log(1/\eps)))^{-1}\int_{\eps}^x \frac{1}{t\log(1/t)} dt, & \qquad \mbox{if { }} \eps < x \leq e^{-1},\\
1, & \qquad \mbox{if { }} e^{-1} < x.
\end{array}
\right .
\ee

Then

\be j'(x) =
(x\log(1/x))^{-1}(\log(\log(1/\eps)))^{-1}1_{[\eps,e^{-1}]}(x)\ee

in the distributional sense. Let $\chi (z) = j(d_U(z))$. Since
$d_U(z)$ is Lipshitz, it has distributional derivatives of norm at
most 1, and thus

\be \lll{nice}|\bar{\dd}\chi(z)|\leq \Big(d_U (z)\log(1/d_U
(z))\log(\log(1/\eps))\Big)^{-1} . \ee

Note also that $0 \leq \chi \leq 1$, with $\chi(w) = 1$ and $\chi(z)
= 0$ whenever $d_U(z)< \eps$. Let $\hat{V}(z) = V(z)\chi(z)$.
Since $H(\dd U,\dd W) \leq \eps$ we see that $\hat{V}(z)$ is 0 on $W^c$ and 1 at $w$, so it gives an upper
bound for $\la_W(w)$. We obtain

\be \bar{\dd}\hat{V} = (\bar{\dd}V)\chi + V(\bar{\dd}\chi) \ee

so that

\be \lll{cool} |\bar{\dd}\hat{V}| \leq |(\bar{\dd}V)| +
|V(\bar{\dd}\chi)| . \ee

We have the estimate $|V(z)| \leq |C d_U(z) \log(1/d_U(z))|$ from
Theorem 7 of Chapter 3 in \cite{fnbook}, where $C$ here depends on $d(K,\dd U)$ and
the diameter of $U$. In light of this and (\ref{nice}), we see that
(\ref{cool}) is bounded above by

\be ||V||_\ff + C(\log(\log(1/\eps)))^{-1} . \ee

Thus, $\la_W (w) \leq \la_U (w) + C(\log(\log(1/\eps)))^{-1}$ for all
$w \in K$. Interchanging the roles of $U$ and $W$ in the above
argument gives the reverse inequality, and shows that $|\la_U (w) -
\la_W (w)| \leq C (\log(\log(1/\eps)))^{-1}$. \qed

\section{Rates of convergence for the hyperbolic density.}

As mentioned in Section \ref{tm}, the Teichm\"{u}ller and hyperbolic densities coincide, up to a constant, on simply connected domains. The following is therefore a corollary to Theorem \ref{big}.

\begin{cor} Suppose that $U$ is a simply connected, bounded domain, and that $K$ is a compact
subset of $U$. Then there is a constant $C$ such that if $W$ is a
simply connected domain with $H(\dd U,\dd W) \leq \eps$ then
$|\rho_U (z) - \rho_W (z)| \leq C (\log(\log(1/\eps)))^{-1}$ for all
$z \in K$, provided that $\eps$ is sufficiently small. The constant $C$ depends on $d(K,\dd U)$ and the diameter of $U$.
\end{cor}

{\bf Remark:} The boundedness requirement may be relaxed in certain cases. For instance, if $U$ is unbounded but $U^c$ contains an open set containing a point $p$, we may apply a M\"{o}bius inversion mapping $p$ to $\ff$. The image under this map is a bounded domain, and the spherical metric is preserved by the inversion up to a constant and is equivalent to the Euclidean metric in the bounded image. Thus, the corollary may be applied, with $H(\dd U,\dd W)$ now being measured in the spherical metric.

\vski

It was also mentioned in Section \ref{tm} that the Teichm\"{u}ller and hyperbolic densities coincide on the largest possible hyperbolic domain, the thrice punctured plane $\hat{\CC} \backslash \{a,b,c\}$. Thus, we may obtain a similar result if $U_n$ and $U$ are thrice punctured planes. However, as the hyperbolic density of a thrice punctured plane is explicitly computable we may obtain a far better rate of convergence when $U_n \lar U$ in boundary, namely $H \log(\frac{1}{H})$, where $H=H(\dd U_n,\dd U)$. Let the hyperbolic metric on $\hat{\CC}\backslash \{a,b,c\}$ be denoted $\rho_{a,b,c}(z)$.

\bt \lll{jets}
Let $U=\CC\backslash \{a,b,c\}$ and let $K$ be a compact set in $U$. Then there is a constant $C$ such that $|\rho_{a,b,c}(z) - \rho_{a',b',c'}(z)|<C\eps \log(\frac{1}{\eps})$ for all $z \in K$ whenever $H(\{a,b,c\}, \{a',b',c'\}) < \eps$ and $\eps$ is sufficiently small. $C$ depends on $\sup_{w \in K} \{|w-a|,|w-b|,|w-c|\}, d(K,\{a,b,c\}),|a-b|, |b-c|,$ and $|a-c|$.
\et

{\bf Proof:}  With no loss of generality we may assume that $|a-a'|,|b-b'|, |c-c'| \leq \eps.$  The triangle inequality

\bea \label{}
&& \nn |\rho_{a,b,c}(z)-\rho_{a',b',c'}(z)| \leq |\rho_{a,b,c}(z)-\rho_{a,b,c'}(z)|
\\ \nn && \hspace{1cm} +|\rho_{a,b,c'}(z)-\rho_{a,b',c'}(z)|+|\rho_{a,b',c'}(z)-\rho_{a',b',c'}(z)|
\eea

implies that we may assume that $a=a'$ and $b=b'.$ \cite{fnbook} contains a proof of the following formula:

\be \lll{ref}
\frac{1}{\rho_{a,b,c}(z)} = \frac{1}{\pi} \int \! \! \int_{\CC}
\frac{|(z-a)(z-b)(z-c)|}{|(w-a)(w-b)(w-c)(w-z)|}
dA(w).
\ee




Using this formula together with the corresponding expression for $\rho_{a,b,c'}(z)$ we have

\bea \label{wes}
&& |\rho_{a,b,c}(z) - \rho_{a,b,c'}(z)| =
\frac{1}{\pi}\rho_{a,b,c}(z) \rho_{a,b,c'}(z)
\\ \nn && \hspace{1cm}  \times \int \! \! \int_{\CC}  \frac{|(z-a)(z-b)|}{|(w-a)(w-b)(w-z)|} \Big(
\frac{|z-c|}{|w-c|} -\frac{|z-c'|}{|w-c'|}
\Big)dA(w).
\eea


Let $M = \sup_{w \in K} \{|w-a|,|w-b|,|w-c|\}$, $d = d(K,\{a,b,c\})$, and $m=\min\{|a-b|, |b-c|,|a-c|, d\}$. Now
\be \label{wel}
\frac{|z-c|}{|w-c|} -\frac{|z-c'|}{|w-c'|}\leq \frac{|z-c|}{|w-c|} -\frac{|z-c'|}{|w-c|}+\frac{|z-c'|}{|w-c|} -\frac{|z-c'|}{|w-c'|}.
\ee

Thus,

\be \label{dir}
|\rho_{a,b,c}(z) - \rho_{a,b,c'}(z)| \leq \rho_{a,b,c}(z) \rho_{a,b,c'}(z)
(I + II)
\ee

where

\be I=\frac{1}{\pi}\int \! \! \int_{\CC}  \frac{|(z-a)(z-b)(c-c')|}{|(w-a)(w-b)(w-c)(w-z)|}dA(w)\ee

and

\be II=\frac{1}{\pi}\int \! \! \int_{\CC}  \frac{|(z-a)(z-b)(z-c')(c-c')|}{|(w-a)(w-b)(w-c)(w-c')(w-z)|}dA(w) .\ee

Note that \be \rho_{a,b,c}(z) \leq \rho_{D(z,d)}(z)=\frac{1}{d}\ee
and similarly  \be\rho_{a,b,c'}(z) \leq \frac{1}{d-\eps}.\ee
Furthermore,
\be I=\frac{1}{\pi}\int \! \! \int_{\CC}  \frac{|(z-a)(z-b)(z-c)|}{|(w-a)(w-b)(w-c)(w-z)|}\frac{|c-c'|}{|z-c|}dA(w) \leq \frac{\eps}{d}\frac{1}{\rho_{a,b,c}(z)}.\ee

We will estimate II by considering two different regions. Let $A=\{w\in \CC: |w-c'|>\frac{d}{2}\}$ and let $B$ be the complement of $A.$ Then

\be II=\frac{1}{\pi}\int \! \! \int_{A}  \frac{|(z-a)(z-b)(z-c')(c-c')|}{|(w-a)(w-b)(w-c)(w-c')(w-z)|}dA(w)+\ee
$$\frac{1}{\pi}\int \! \! \int_{B}  \frac{|(z-a)(z-b)(z-c')(c-c')|}{|(w-a)(w-b)(w-c)(w-c')(w-z)|}dA(w).$$

We have

\be \frac{1}{\pi}\int \! \! \int_{A}  \frac{|(z-a)(z-b)(z-c')(c-c')|}{|(w-a)(w-b)(w-c)(w-c')(w-z)|}dA(w)=\ee
$$\frac{1}{\pi}|c-c'|\int \! \! \int_{A}  \frac{|(z-a)(z-b)(z-c)|}{|(w-a)(w-b)(w-c)(w-z)|}\frac{|z-c'|}{|w-c'||z-c|}dA(w)$$ $$\leq \frac{2\eps (M+\eps)}{d^2\rho_{a,b,c}(z)},$$

and for the integral over $B$ we obtain

\be \frac{1}{\pi}\int \! \! \int_{B}  \frac{|(z-a)(z-b)(z-c')(c-c')|}{|(w-a)(w-b)(w-c)(w-c')(w-z)|}dA(w)=\ee
$$\frac{1}{\pi}\int \! \! \int_{B}  \frac{|(c'-a)(c'-b)(c'-c)|}{|(w-a)(w-b)(w-c)(w-c')|}\frac{|(z-a)(z-b)(z-c')|}{|(c'-a)(c'-b)(w-z)|}dA(w)$$ $$\leq \frac{M^2 (M+\eps)}{(\frac{d}{2}-\eps)(m-\eps)^2\rho_{a,b,c}(c')}$$

where we have used

\be \lll{ref1}
\frac{1}{\rho_{a,b,c}(c')} = \frac{1}{\pi} \int \! \! \int_{\CC}
\frac{|(c'-a)(c'-b)(c'-c)|}{|(w-a)(w-b)(w-c)(w-c')|}
dA(w) .
\ee

Now let $f$ be the Mobius transformation

\be f(z)=\frac{(z-c)(a-b)}{(z-b)(a-c)}.\ee
We have $f(c)=0, f(b)=\infty$ and $f(a)=1.$  Thus
\be \rho_{a,b,c}(z)=\rho_{0,1}(f(z))|f'(z)|=\rho_{0,1}(f(z))\frac{|(a-b)(b-c)|}{|(a-c)(z-b)^2|}\ee
where we are using the notation $\rho_{0,1}(z)$ as a shorthand for $\rho_{0,1,\ff}(z)$. For $z \in W,$
\be \frac{|(a-b)(b-c)|}{|(a-c)(z-b)^2|}\geq \frac{m^2}{2Md^2}. \ee

Note that $|f(z)| = \frac{|(z-a)(a-b)|}{|(z-b)(a-c)|} \leq \frac{2M^2}{dm}$.
Theorem 14.3.1 in \cite{nlbook} shows then that

\begin{equation} \label{}
\rho_{0,1}(f(z)) \geq \frac{1}{2|\frac{2M^2}{dm}|\log|\frac{2M^2}{dm}| + 10 \frac{2M^2}{dm}}.
\end{equation}

Thus, $\rho_{a,b,c}(z)$ is bounded below by a constant depending on $m,M,d$, and we conclude that
\be \frac{1}{\pi}\int \! \! \int_{A}  \frac{|(z-a)(z-b)(z-c')(c-c')|}{|(w-a)(w-b)(w-c)(w-c')(w-z)|}dA(w)\leq C\eps.\ee
Similarly,
\be \rho_{a,b,c}(c')=\rho_{0,1}(f(c'))|f'(c')|=\rho_{0,1}(f(z))\frac{|(a-b)(b-c)|}{|(a-c)(c'-b)^2|}.\ee
We have
\be \frac{|(a-b)(b-c)|}{|(a-c)(c'-b^2)|}\geq \frac{m^2}{2M(2M+\eps)^2}\ee and
\be f(c')=\frac{|(a-b)(c'-c)|}{|(a-c)(c'-b)|}\leq \frac{2M\eps}{m(m-\eps)}<\frac{1}{2}\ee for sufficiently small $\eps.$   Corollary 14.4.1 in \cite{nlbook} implies
\be \frac{1}{\rho_{01}(f(c'))} \leq 17 |f(c')|\log (\frac{1}{|f(c')|}) \leq C \eps \log (\frac{1}{\eps}).\ee
Thus, we have
\be \frac{1}{\pi}\int \! \! \int_{B}  \frac{|(z-a)(z-b)(z-c')(c-c')|}{|(w-a)(w-b)(w-c)(w-c')(w-z)|}dA(w)\leq C\eps\log (\frac{1}{\eps}),\ee
 and this proves the theorem. \qed

In the prior theorem we have not allowed any of $a,b,c$ to be $\ff$. If we allow one of the points to be $\ff$, using the notation $\rho_{a,b}$ in place of $\rho_{a,b,\ff}$, we can obtain the following.

\begin{theorem} \label{suck}
Let $U = \CC \backslash \{a,b\}$ and let $K$ be a compact set in $U$. Then there is a constant $C$ such that $|\rho_{a,b}(z)-\rho_{a',b'}(z)| < C \eps \log(\frac{1}{\eps})$ for all $z \in K$ whenever $H(\{a,b\},\{a',b'\})<\epsilon$. The constant $C$ depends on  $\sup_{w \in K} \{|w-a|,|w-b|\}$, $d(K,\{a,b\})$ and $|a-b|$.
\end{theorem}

{\bf Proof:} Let $\gamma$ be a line segment joining $K$ and $\{a,b\}$ with the length of $\gamma$ equal to $d(K,\{a,b\})$. Let $r$ be the midpoint of $\gamma$. For any point $z$ in $K$ the M\"{o}bius transformation $g(z)=\frac{1}{(z-r)}$ satisfies

\begin{equation} \label{}
\rho_{a,b}(z)=\rho_{1/(a-r),1/(b-r),0}(1/(z-r))\frac{1}{|z-r|^2}
\end{equation}

and

\begin{equation} \label{}
\rho_{a',b'}(z)=\rho_{1/(a'-r),1/(b'-r),0}(1/(z-r))\frac{1}{|z-r|^2}.
\end{equation}

Theorem \ref{suck} now follows from Theorem \ref{jets} applied to the compact set $\{\frac{1}{(z-r)}: z \in K\}$, using the the point $\frac{1}{(z-r)}$ in place of $z$. \qed

In the case of the the unit disc, we can obtain a $H(\dd U_n,\dd \DD)$ rate of convergence.

\bt
If $K$ is a compact subset of $\DD$, then there is a constant $C$ depending on $K$ such that $|\rho_{\DD}(z)-\rho_W(z)| \leq C\eps$ for all $z$ in $K$ whenever $W$ contains $K$ and $H(\dd \DD, \dd W) \leq \eps$ for $\eps$ sufficiently small.
\et

{\bf Proof:} For $r>0$ let $r\DD = \{|z|<r\}$. The hyperbolic density on $r\DD$ is well known to be given by $\rho_{r\DD}(z) = \frac{r}{r^2-|z|^2}$ (see \cite{fnbook}). Thus,

\be \lll{}
\begin{split}
|\rho_{r\DD}(z)-\rho_{\DD}(z)| & = \frac{|r(1-|z|^2) - (r^2-|z|^2)|}{(1-|z|^2)(r^2-|z|^2)}
\\ & \leq \frac{(r|1-r|+ |z|^2|1-r|)}{(\min (1,r^2) - |z|^2)^2}\leq C|1-r|
\end{split}
\ee

for $|z|$ uniformly bounded below $r$. The result now follows by the monotonicity property of the hyperbolic density, as when $H(\dd \DD, \dd W) \leq \eps$ we must have $(1-\eps)\DD \subseteq W \subseteq (1+\eps)\DD$. \qed

{\bf Remark:} The examples $\{|z|<1+\eps\}$ and $\{|z|<1-\eps\}$ show that this rate of convergence can not be improved.

\vski

Though not directly related to the results given in this section, we would be remiss if we did not mention the one prior result we have seen concerning the rate of convergence of the hyperbolic density. In \cite{beard} the rate of convergence of $\rho_U(x)$ was determined for $x \in \reals$ and $U$ of the form $\{|x| < l, |y| < \frac{\pi}{2}\}$ for changing values of $l$.

\section{Remarks on the three-point density.}

Suppose $U$ is an arbitrary domain in $\CC$. If $a,b,c \in U^{c}$, then $\rho_{a,b,c}(z) \leq \rho_U(z)$ for all
$z \in U$ by the monotonicity of the hyperbolic density. We can
define a new density on $U$ by setting

\be h_U(z) = \sup_{a,b,c \in U^{c}} \rho_{a,b,c}(z) . \lll{8} \ee

This was first done in \cite{fnart}, and we shall refer to this quantity as the {\it three-point density}. We clearly have $h_U \leq \rho_U$, and it is true, though less
obvious, that equality holds only when $U$ is itself the thrice
punctured sphere. It was shown initially in \cite{fnart} that
there is a positive universal constant $C$ such that

\be \label{dars}
\rho_U \leq Ch_U,
\ee

so that the two densities are equivalent. In \cite{minma}, it was shown that $h_U$ is a continuous density which is M\"{o}bius invariant. \cite{suvu} gave an explicit constant for \rrr{dars}, and \cite{bets} worked to improve the constant and also calculated $h_\DD$. In relation to the convergence of densities, we have the following proposition.

\begin{prop} If $z \in U$ and $U_n \lar U$ in boundary, then $h_{U_n}(z) \lar h_U(z)$.
\end{prop}

Proof: If $a,b,c \in U^c$ we can choose $a_n,b_n,c_n \in U_n^c$
which converge to $a,b,c$ respectively. By Lemma \ref{old},

\be \rho_{a_n,b_n,c_n}(z) \lar \rho_{a,b,c}(z) . \lll{9} \ee

It follows from this that

\be \underline{\lim}h_{U_n}(z) \geq h_U(z) . \lll{10} \ee

For the reverse inequality, we can choose $a_n,b_n,c_n \in U_n^c$
such that

\be \overline{\lim} \rho_{a_n,b_n,c_n}(z) = \overline{\lim}
h_{U_n}(z) . \lll{10a} \ee

After passing to subsequences several times if necessary we may
assume $a_n,b_n,c_n \lar a,b,c$ respectively. Since $U_n \lar U$ we
can choose $a_n',b_n',c_n' \in U^c$ close to $a_n,b_n,c_n$, such
that $a_n',b_n',c_n' \lar a,b,c$. $U^c$ is closed, so $a,b,c \in U^c$.
If $a,b,c$ are distinct, then applying Lemma \ref{old} we see that

\be \overline{\lim}h_{U_n}(z) \leq h_U \lll{11}(z) \ee

completing the proof. It remains only to see that $a,b,c$ must be all distinct, since if not then we would approach a pole of order two in \rrr{ref} as $n \lar \ff$. This would force $\rho_{a_n,b_n,c_n}(z) \lar
0$, contradicting (\ref{10}) and \rrr{10a}. \qed

It would seem that Theorem \ref{jets} was ideally suited for deducing a rate of convergence result for $h_{U_n}$. In fact, it was shown in \cite{minma} that the supremum in \rrr{8} is always attained for some triple $a,b,c$ in $\dd U$, so if $U$ is bounded we may assume $a,b,c$ are bounded in \rrr{8} as well. The difficulty, however, lies in the fact that the constant in Theorem \ref{jets} depends in part on $(\min\{|a-b|, |b-c|,|a-c|\})^{-1}$. If a point is fixed in $U$ then there are $a,b,c \in \dd U$ such that $h_U(z) = \sup_{a,b,c \in U^{c}} \rho_{a,b,c}(z)$. We may then apply Theorem \ref{jets} to obtain

\be \label{}
h_{U_n}(z) \geq h_U(z) - C \sqrt{H(\dd U, \dd U^n)}
\ee

for $H(\dd U, \dd U^n)$ sufficiently small. However, different points in a compact set $K$ determine different optimal triples of points $a,b,c$, and we do not currently have a way to bound $\min\{|a-b|, |b-c|,|a-c|\}$ from below. For the argument to show $\limsup h_{U_n}(z) \leq h_U(z)$, we know that $a_n,b_n,c_n$ converge, but we do not know how far the limit points are from each other. Obtaining a theorem on the rate of convergence of $h_U$ would seem therefore to necessitate understanding how the optimal points $a,b,c$ are situated in the plane for given $z$ and $U$.

\section{Further questions}

It may be of interest for applications to explicitly calculate the constants in the results given above. It would also be interesting to know whether the rates of convergence given are the best possible. Except where stated, we do not know whether this is the case. Perhaps there are results similar to the ones in this paper for any of a number of other densities, for instance the Carath\'{e}odory density or Kobayashi density in higher dimensions. Of course, finding an analog of Corollary 1 for domains which are not simply connected would be desirable as well.

\section{Acknowledgements}

We would like to thank Fred Gardiner for many helpful conversations, as well as an anonymous referee for a careful reading. 


\def\noopsort#1{} \def\printfirst#1#2{#1} \def\singleletter#1{#1}
   \def\switchargs#1#2{#2#1} \def\bibsameauth{\leavevmode\vrule height .1ex
   depth 0pt width 2.3em\relax\,}
\makeatletter \renewcommand{\@biblabel}[1]{\hfill#1.}\makeatother

\bibliographystyle{alpha}
\bibliography{CAbib}

{\it Nikola Lakic 

nikola.lakic@lehman.cuny.edu 

Lehman College, CUNY 

Bronx, NY 10468 U.S.A }

\vski

{\it Greg Markowsky

gmarkowsky@gmail.com

Monash University

Victoria 3800 Australia
\end{document}